\documentclass[11pt,a4paper]{article}

\usepackage[english]{babel}
\usepackage{amsmath,amssymb,amsfonts,a4}
\usepackage{eucal,bbm,stmaryrd}
\usepackage{hyperref}
\usepackage{xcolor,tikz}

\begin{document}

\numberwithin{equation}{section}

\newtheorem{Thm}{Theorem}
\newtheorem{Prop}{Proposition}
\newtheorem{Def}{Definition}
\newtheorem{Lem}{Lemma}
\newtheorem{Rem}{Remark}
\newtheorem{Cor}{Corollary}
\newtheorem{Con}{Conjecture}

\newcommand{\Thmautorefname}{Theorem}
\newcommand{\Propautorefname}{Proposition}
\newcommand{\Defautorefname}{Definition}
\newcommand{\Lemautorefname}{Lemma}
\newcommand{\Remautorefname}{Remark}
\newcommand{\Corautorefname}{Corollary}
\newcommand{\Conautorefname}{Conjecture}

\newcommand{\Pf}{\noindent{\bf Proof: }}
\newcommand{\qed}{\hspace*{3em} \hfill{$\square$}}

\newcommand{\N}{\mathbbm{N}}
\newcommand{\Z}{\mathbbm{Z}}

\newcommand{\E}{\mathbbm{E}}
\renewcommand{\P}{\mathbbm{P}}

\newcommand{\Bin}{\mathit{Bin}}

\newcommand{\la}{\lambda}
\newcommand{\La}{\Lambda}
\newcommand{\si}{\sigma}
\newcommand{\al}{\alpha}
\newcommand{\be}{\beta}
\newcommand{\ep}{\epsilon}
\newcommand{\ga}{\gamma}
\newcommand{\de}{\delta}
\newcommand{\De}{\Delta}
\newcommand{\ph}{\varphi}
\newcommand{\om}{\omega}
\renewcommand{\th}{\theta}
\newcommand{\vth}{\vartheta}
\newcommand{\ka}{\kappa}

\newcommand{\cZ}{\mathcal{Z}}
\newcommand{\cT}{\mathcal{T}}
\newcommand{\cX}{\mathcal{X}}
\newcommand{\cY}{\mathcal{Y}}
\newcommand{\cW}{\mathcal{W}}
\newcommand{\tcW}{\tilde{\cW}}
\newcommand{\cV}{\mathcal{V}}
\newcommand{\cC}{\mathcal{C}}
\newcommand{\cN}{\mathcal{N}}
\newcommand{\da}{\dagger}
\newcommand{\Mda}{M^{\dagger}}
\newcommand{\tM}{\tilde{M}}
\newcommand{\xda}{x_\da}
\newcommand{\xst}{x_*}

\newcommand{\bz}{\bar{z}}
\newcommand{\Aut}{Aut}

\newcommand{\Nm}{N_{\wedge}}

\newcommand{\stm}{\setminus}
\newcommand{\lra}{\leftrightarrow}
\newcommand{\lraa}[1]{\stackrel{#1}{\lra}}
\newcommand{\Ra}{\Rightarrow}
\newcommand{\Lra}{\Leftrightarrow}

\newcommand{\sq}{ \, \square \,}

\makeatletter 
\newcommand{\nobr}{\par\nobreak\@afterheading\vspace{0.3 cm}} 
\makeatother

\thispagestyle{plain}
\title{Monotonicity properties\\ 
for Bernoulli percolation on layered graphs\\ 
- a Markov chain approach}
\author{Philipp König and Thomas Richthammer
\footnote{Institut f\"ur Mathematik, Universit\"at Paderborn}}

\maketitle
\begin{abstract}
A layered graph $G^\times$ is the Cartesian product 
of a graph $G = (V,E)$ with the linear graph $\Z$, 
e.g. $\Z^\times$ is the 2D square lattice $\Z^2$. 
For Bernoulli percolation with parameter $p \in [0,1]$ 
on $G^\times$ one intuitively would expect that 
$\P_p((o,0) \lra (v,n)) \ge \P_p((o,0) \lra (v,n+1))$
for all $o,v \in V$ and $n \ge 0$. 
This is reminiscent of the better known bunkbed conjecture. 
Here we introduce an approach to the above monotonicity conjecture 
that makes use of a Markov chain building the percolation pattern layer by layer. 
In case of finite $G$ we thus can show 
that for some $N \ge 0$ the above holds for all $n \ge N$ $o,v \in V$ and $p \in [0,1]$.  
One might hope that this Markov chain approach could be useful 
for other problems concerning Bernoulli percolation on layered graphs. 

\end{abstract}

%



\section{Introduction}

Bernoulli (bond) percolation is a stochastic process that 
has been studied extensively since it was introduced in the 1950s. 
For a given connected (simple) graph $G = (V(G),E(G))$ 
and a parameter value $p \in (0,1)$ 
we consider independent random variables $\cZ_e, e \in E(G)$, 
such that $\P(\cZ_e = 1) = p$ and $\P(\cZ_e = 0) = 1-p$. 
Often the value $p$ is indicated in the notation for probabilities in the form of $\P_p$. 
$e$ is called open if $\cZ_e = 1$, and $e$ is called closed if 
$\cZ_e = 0$. 
The family of random variables $(\cZ_e)_{e \in E(G)}$ determines a random 
subgraph $G_\cZ$ of $G$ via $V(G_\cZ) := V$ and $E(G_\cZ) := \{e \in E: \cZ_e = 1\}$. 
Usually $G$ is assumed to be regular in some sense (e.g. a lattice). 
For variants of this process and various results and open questions in Bernoulli percolation see the standard reference \cite{G}. 
For $u,v \in V(G)$ we are interested in the event $\{u \lra v\}$
that $u,v$ are in the same connected component of $G_\cZ$. 
Interpreting Bernoulli percolation as a model for the spread of a disease, 
where an open bond transmits the disease and a closed bond does not,  
$\{o \lra v\}$ represents the event 
that $v$ is infected by a disease originating in some fixed vertex $o \in V(G)$. 
It is not unreasonable to expect that vertices closer to $o$ 
(in some sense that has to be specified) 
are more likely to be infected than vertices further away, i.e. 
we have some sort of spatial monotonicity property of the connectivity function 
(or two-point function) $f_p(v) := \P_p(o \lra v)$. 
In the following we present three settings, where 
the above monotonicity property can be formulated in a precise way.  

\begin{Con} \label{Con:Zd}
Monotonicity of the connectivity function in case of $\Z^d$. 
Let $d \ge 1$ and $p \in (0,1)$. 
On the vertex set of the graph $\Z^d$
(in which any two vertices that have a Euclidean distance of $1$
are connected by an edge)
we consider the componentwise (partial) order $\le$ (i.e. 
$u \le v$ iff $u_i \le v_i$ for all $i \in \{1,...,d\}$). 
For Bernoulli percolation on $\Z^d$ with parameter $p$ we have 
\begin{equation} \label{equ:Zd}
\forall o,u,v \in \Z^d: \quad 
o \le u \le v \quad \Ra \quad \P_p(o \lra u) \ge \P_p(o \lra v). 
\end{equation}
\end{Con}
This conjecture is mathematical folklore, and we don't know its precise origin. 
Similar natural monotonictiy conjectures can be made for other 
regular graphs, and there are variants of the above conjecture for $\Z^d$:  
e.g. one might want to consider edge probabilities $p_i, i \in \{1,...,d\}$, 
depending on the direction of the edges. 
In case of $d = 1$ the above conjecture is trivially true, 
but already the case $d = 2$ is open (even if it seems to be 
very compelling to assume it is true). 
The only pertaining result we know of is the result of \cite{LPS}, where the above monotonicity is established for  $o,u,v \in \Z \times \{0\}^{d-1}$ 
provided that $p$ is sufficiently close to $0$. 

A classical generalization of the above conjecture is the so called 
bunkbed conjecture that applies to graphs of the following structure: 

\begin{Def} \label{Def:bunkbed}
Let $G$ be a graph and $T \subset V(G)$. 
The {\bf bunkbed graph} $G_T$ 
is given by $V(G_T) = V(G) \times \{0,1\}$ and $E(G_T) = E^h(G_T) \cup E^v(G_T)$, where 
$E^h(G_T) = \{(u,l)(v,l): uv \in E(G), l \in \{0,1\}\}$ 
is the set of horizontal edges 
and $E^v(G_T) = \{(u,0)(u,1): u \in T\}$ is the set of vertical edges.  
\end{Def}
Loosely speaking a bunkbed graph consists of two copies of $G$ with 
vertices in the two copies connected to each other (by a vertical edge) iff they are copies of the same vertex in $T$.
The name 'bunkbed graph' comes from thinking of the two copies of $G$ as the two beds of a bunkbed and the vertical edges as the poles connecting the two beds (where the pole positions are given by $T$). 
 
\begin{Con} \label{Con:bunkbed}
Bunkbed conjecture. Let $G = (V,E)$ be a finite connected graph, 
$T \subset V$, and let $p \in (0,1)$. 
For Bernoulli percolation on the bunkbed graph $G_T$ with parameter $p$ 
we have 
\begin{equation}\label{equ:bunkbed}
\forall o,v \in V: \; 
\P_p((o,0) \lra (v,0)) \ge \P_p((o,0) \lra (v,1)). 
\end{equation}
\end{Con}
This is a natural extension of \autoref{Con:Zd} and can be traced back to P.W. Kasteleyn (1985) (as noted in \cite{BK}). 
The bunkbed conjecture also appears in several variants. 
Sometimes only the case $T = V$ is considered. 
More general, one may want to consider the case of infinite connected graphs $G$, or the case of nonconstant edge probabilities, 
say $p_e \in [0,1]$, $e \in E$ for the horizontal edges
and $p_u \in [0,1]$, $u \in V$ for the vertical edges 
(and indeed, often only the case of deterministic vertical edges, i.e. $p_u \in \{0,1\}$, is considered.) 
Again it is very compelling to assume the conjecture is true, 
but while it has received a lot of attention, 
only few partial results have been obtained, see  
\cite{Ri2} and the references therein. 
In this paper we we would like to focus on the setting 
of layered graphs
that is somewhat in between the two settings considered above. 

\begin{Def} \label{Def:layer}
Let $G = (V,E)$ be a graph. 
The {\bf layered graph} $G^\times$  
is the Cartesian product of $G$ and $\Z$, 
i.e. $V(G^\times) = G \times \Z$ and $E(G^\times) = E^h(G^\times) \cup E^v(G^\times)$, 
where $E^h(G^\times) = \{(u,l)(v,l) \!: uv \in E, l \in  \Z\}$ is 
the set of horizontal edges  
and $E^v(G^\times) = \{(u,l-1)(u,l) \!: u \in V, l \in \Z\}$ is the 
set of vertical edges. 
\end{Def}
Loosely speaking the layered graph $G^\times$ consists of biinfinitely many copies of $G$ 
(forming the layers of $G^\times$) 
with vertices in two distinct layers connected to each other 
(by a vertical edge) iff they are copies of the same vertex in $V$ and 
the two layers are adjacent. 
To illustrate this definition we note that $\Z^d$ can be viewed 
as the layered graph $(\Z^{d-1})^\times$.

\begin{Con} \label{Con:layer}
Monotonicity of the connectivity function in case of layered graphs. 
Let $G = (V,E)$ be a finite connected graph, and $p \in (0,1)$. 
For Bernoulli percolation on the layered graph $G^\times$ with parameter $p$ we have 
\begin{equation} \label{equ:layer}
\forall o,v \in V \; \forall 0 \le n \le n': 
\;
\P_p((o,0) \lra (v,n)) \ge \P_p((o,0) \lra (v,n')). 
\end{equation}
\end{Con}
This conjecture seems to be new. 
Again, more generally one may want to consider the case of 
infinite connected graphs $G$, or the case of 
nonconstant edge probabilities, 
say $p_e \in [0,1]$, $e \in E$ for the horizontal edges (in each layer)
and $p_u \in [0,1]$, $u \in V$ for the vertical edges. 
As we will show, the above conjecture is somewhat in between 
the two preceding conjectures in terms of generality. 
We will use the layered structure to investigate 
weaker and stronger versions of the above conjecture 
in this paper, namely the monotonicity of the expected number of 
infected vertices per layer and the monotonicity of the probability 
for an infection pattern within a layer. 
While our results concerning the former rely on combinatorial arguments and classical tools from percolation theory, 
our results concerning the latter rely on a combination of classical tools with Markov chain methods: 
Building the percolation structure layer by layer, we obtain a Markov chain, 
for which asymptotics for transition probabilities can be used to obtain monotonicity results. 
As far as we know this is a novel approach to investigating 
percolation problems in general, 
and we hope that this new perspective will be of some use in the future. 

The paper is organized as follows. 
In Section \ref{Sec:results} we present our results 
and explain how they are connected to the above conjectures. 
In Section \ref{Sec:comparison} we prove statements on the comparison of monotonicity properties. 
In Section \ref{Sec:Markov} we investigate the Markov chain 
on infection patterns. 
In Sections \ref{Sec:patterns} we prove our main result on the monotonicity of the probability for an infection pattern within a layer.
In Section \ref{Sec:Number} we prove our result  
on the monotonicity of the expected number of infected points per layer.

\section{Results} \label{Sec:results}

As a motivation for our interest in Conjecture \ref{Con:layer} 
let us note its relation to Conjectures \ref{Con:Zd} and \ref{Con:bunkbed}. Recall that the Cartesian product of graphs 
$G_1,...,G_d$ with $G_i  = (V_i,E_i)$ is the graph $G = (V,E)$ 
with $V = V_1 \times ... \times V_d$ and 
$E = \{ (v_1,..,v_d)(w_1,...,w_d): \exists i: 
v_iw_i \in E_i, \forall j \neq i: v_i = w_i \in V_i\}$. 
For $d \ge 1$ let $\Z^d$,  $L_k^d$ and $C_k^d$ denote 
the Cartesian products of $d$ copies of $\Z$, the linear graph $L_k$ and the cycle graph $C_k$ on $k \ge 2$ vertices respectively. 
Here we adopt the convention that $C_2 := L_2$ 
in order to avoid multi-edges. 

\begin{Prop} \label{Prop:bunkbed}
\autoref{Con:bunkbed} (on bunkbed graphs) implies \autoref{Con:layer} (on layered graphs), 
and \autoref{Con:layer}  
implies \autoref{Con:Zd} (on $\Z^d$). 
More precisely: 
\begin{itemize}
\item[(a)] Let $p \in (0,1)$. 
If \eqref{equ:bunkbed} holds for Bernoulli percolation 
with parameter $p$ on every bunkbed graph $G_T$ 
(i.e. for every finite connected graph $G$ and every $T \subset V(G)$), 
then \eqref{equ:layer} holds for Bernoulli percolation with parameter $p$ 
on every layered graph $G^\times$ 
(i.e. for every finite connected graph $G$).
\item[(b)] Let $p \in (0,1)$ and $d \ge 2$. 
If \eqref{equ:layer} holds for Bernoulli percolation with parameter 
$p$ on the layered graphs $(C_k^{d-1})^\times$ (i.e. 
for every $k \ge 2$), 
then \eqref{equ:Zd} holds for Bernoulli percolation with parameter $p$ on 
$\Z^d$. (The same is true for $L_k^{d-1}$ instead of $C_k^{d-1}$.) 
\end{itemize}
\end{Prop}
The proof will be given in Section \ref{Sec:comparison}. 
Since $\Z^2$ is the most prominent graph on which Bernoulli percolation is studied, 
in view of (b) we are mostly interested in the layered graphs 
$C_k^\times$. 
While our results are fairly general, 
one might want to keep this particular example in mind. 

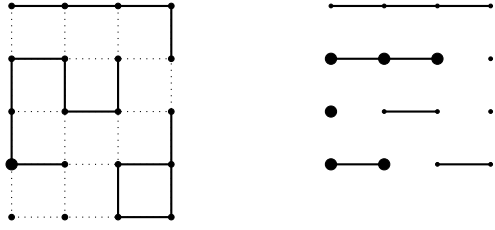
\begin{figure}[htb!] 
\centering
\vspace*{0.5 cm}
\begin{tikzpicture}[scale = 0.7]
\draw[dotted] (0,-1) grid (3,3);
\foreach \x in {0,1,2,3} 
\foreach \y in {-1,0,1,2,3}
\draw[fill] (\x,\y) circle (1.5 pt);
\draw[fill] (0,0) circle (3 pt);
\draw[thick] (1,0)--(0,0)--(0,1)--(0,2)--(1,2)--(1,1)--(2,1)--(2,2);
\draw[thick] (2,-1)--(2,-1)--(2,0)--(3,0)--(3,1);
\draw[thick] (2,-1)--(3,-1)--(3,0);
\draw[thick] (0,3)--(1,3)--(2,3)--(3,3)--(3,2);
\foreach \x in {6,7,8,9}
\foreach \y in {0,1,2,3}
\draw[fill] (\x,\y) circle (1 pt);
\foreach \x in {(6,0),(7,0), (6,1), (6,2),(7,2),(8,2)}
\draw[fill] \x circle (3 pt);
\draw[thick] (6,0)--(7,0) (8,0)--(9,0) (7,1)--(8,1) (6,2)--(7,2)--(8,2) 
(6,3)--(7,3)--(8,3)--(9,3); 
\end{tikzpicture}
\caption{The left hand side shows a realization of Bernoulli percolation on $L_4^\times$
restricted to $\{0,1,2,3\} \times \{-1,0,1,2,3\}$. 
The origin $(0,0)$ is marked with a thick dot, open edges are drawn, closed edges are dotted. 
The right hand side shows the corresponding realization of $\cX_0,...,\cX_3$. Infected vertices 
are marked with a thick dot, connected vertices are joined by lines. With the notation introduced below 
we have $\cX_0 = \{\{*,0,1\},\{2,3\}\}, \cX_1 = \{\{*,0\},\{1,2\},\{3\}\}, \cX_2 = \{\{*,0,1,2\}, \{3\}\}, \cX_3 = \{\{*\},\{0,1,2,3\}\}$. We note that in $\cX_1$ $0$ is not connected to $1$, since here only edges
in $E_{..1}$ may be used in connecting paths. 
}
\label{Fig:Markov}
\end{figure}

From now on we consider the layered graph $G^\times$ 
for some fixed connected graph $G = (V,E)$, 
and we will suppress dependencies on $G$ (and $V$, $E$) 
in our notations. 
We consider a fixed vertex $o \in V$ such that 
$(o,0)$ is considered to be the origin of an infection. 
On $G^\times$ the translation 
$$
\tau: V \times \Z \to V \times \Z, \tau(v,n) := (v,n+1)
$$ 
is a graph automorphism, and we note that Bernoulli-percolation is 
translation invariant in that 
$(\cZ_e)_{e \in E(G^\times)} \sim (\cZ_{\tau(e)})_{e \in E(G^\times)}$
We also introduce notation to refer to vertices and edges 
in certain layers of $G^\times$. For $n \in \Z$ 
let 
$$
V_n := V \times \{n\}, \;  
E_n^h :=  \{(u,n)(v,n) \!: uv \in E\},  \; 
E_n^v := \{(u,n-1)(u,n) \!: u \in V\}
$$ 
and $E_n := E_n^h \cup E_n^v$. 
We think of $V_n$, $E_n^h$, $E_n^v$
as the vertices, the horizontal edges and the vertical edges 
of the $n$-th layer of $G^\times$. 
We also write 
$$
E_{k..n} := \bigcup_{k \le m \le n} E_m, \quad  
E_{..n} := \bigcup_{m \le n} E_m \quad \text{ and } \quad 
E_{n..} := \bigcup_{m \ge n} E_m.  
$$
Our main tool for analyzing connection probabilities 
for Bernoulli percolation on a layered graph 
will be a Markov chain $\cX_n$ that makes use of the layered structure. 
Informally $\cX_n$ describes which vertices in layer $V_n$ are connected to each other 
and which are infected (i.e. connected to $(o,0)$) by means of open paths 
that are fully contained in $E_{..n}$. 
We will formalize this in the following, but it may be more instructive
 to look at the example illustrated in Figure \ref{Fig:Markov}. 

\begin{Def}
Let $* \notin V$ be an abstract symbol. 
A partition $x$ of $V \cup \{*\}$ will be called an (infection) pattern on $V$. 
Let $\sim_x$ denote the corresponding equivalence relation on $V \cup \{*\}$. 
Let $M$ denote the set of all patterns on $V$,  
and let
$$
M^* := \{x \in M: \exists u \in V: u \sim_x *\} \quad \text{ and } \quad  
\Mda := M \setminus M^*. 
$$
\end{Def}
By definition, a pattern $x$ on $V$ 
is the set of equivalence classes of the equivalence relation $\sim_x$ on  $V \cup \{*\}$.  
If for $u,v \in V$ we have 
$u \sim_x v$
we think of $u,v$ as connected w.r.t.~$x$,
and if we have $u \sim_x *$ we think of  $u$ as infected w.r.t.~$x$. 
$M^*$ and $M^\da$ can be interpreted as 
the sets of patterns on $V$
with and without infected vertices respectively. 
In the following definition and later on we will consider 
patterns in certain layers of $G^\times$ that are
built from given (random or nonrandom) edge sets: 
Such an edge set defines a graph and a corresponding connectivity relation on this graph, and the connectivity 
relation restricted to the layer under consideration 
defines a corresponding pattern in this layer, provided 
we specify, where the infection originates. 
We note that a random edge set gives a random connectivity relation and thus a random pattern. 
\begin{Def} We consider percolation on the layered graph $G^\times$. 
Let $n \ge 0$. 
\begin{itemize}
\item[(a)] Let $\lra_{\le n}$ denote the connectivity relation 
induced by percolation on $E_{..n}$.  
Let $\cX_n$ denote the pattern in $V_n$ induced by $\lra_{\le n}$ with the infection originating in $o' := (o,0)$. 
\item[(b)] 
Let $\cW_n \!:= |\{v' \! \in V_n \!\!: v' \! \lra \! o'\}|$ denote the number of infected vertices of $V_n$.
\end{itemize}
\end{Def} 
For clarity we spell out part (a) explicitly: We have 
$u'  \lra_{\le n} v'$ iff there is path from $u'$ to $v'$ 
consisting of open bonds from $E_{..n}$, and we have 
$\cX_n = x$ iff 
\begin{align*}
\forall u,v \in V \!: u \sim_x v  \Lra (u,n) \lra_{\le n} (v,n)
\; \text{ and } \; 
\forall u \in V: u \sim_x * \Lra  (u,n) \lra_{\le n} o'. 
\end{align*}
We note that $\{\cX_n \in M^*\} = \{o'\lra_{\le n} V_n\} = \{o' \lra V_n \} = \{\cW_n > 0 \}$. 
We now relate \autoref{Con:layer} to monotonicity 
properties of the above objects. 

\begin{Prop} \label{Prop:monodiff}
We consider percolation on $G^\times$ for some finite connected graph $G$. 
For all $n \ge 0$ and $p \in (0,1)$ 
we have the implications (i) $\Ra$ (ii) $\Ra$ (iii) 
of the following statements: 
\begin{enumerate}
\item[(i)] For all $x \in M^*$ we have 
$\P_p(\cX_n = x) \ge \P_p(\cX_{n+1} = x)$. 
\item[(ii)] 
For all $v \in V$ we have $\P_p((o,0) \lra (v,n)) \ge 
\P_p((o,0) \lra (v,n+1))$.
\item[(iii)]
We have $\E_p(\cW_n) \ge \E_p(\cW_{n+1})$.  
\end{enumerate}
\end{Prop}
This will be proved in Section \ref{Sec:comparison}. 
The advantage of working with infection patterns instead of the connection events we are mainly interested in, is the following: 

\begin{Prop} \label{Prop:Markov}
We consider percolation on $G^\times$ for some finite connected graph $G$ 
with parameter $p \in (0,1)$. 
$\cX_n, n \ge 0$, is a (time-homogeneous) Markov chain w.r.t. $\P_p$. 
\end{Prop}
This will be shown in Section \ref{Sec:Markov}. 
An obvious disadvantage of working with infection patterns is 
that monotonicity property $(i)$ cannot be expected 
to hold  for all $n \ge 0$. 
Indeed, for any pattern $x \in M^*$ such that $o \not \sim_x *$ 
we have $\P_p(\cX_0 = x) = 0$, whereas for many such patterns 
we have $\P_p(\cX_n = x) > 0$ for some $n > 0$. 
However, it turns out that it is possible to prove 
something slightly weaker, which is our main result:  
\begin{Thm} \label{Thm:pattern}
For any finite connected graph $G$ there is an $N \ge 0$ such that 
for percolation on the layered graph $G^\times$ we have 
\begin{equation} \label{equ:monoNp}
\forall n \ge N \; \forall p \in (0,1) \; \forall x \in M^*: 
\P_p(\cX_n = x) \ge \P_p(\cX_{n+1} = x).
\end{equation}
In particular this implies that for all $n \ge N$, $p \in [0,1]$ and $v \in V$ we have 
$$
\P_p((o,0) \lra (v,n)) \ge 
\P_p((o,0) \lra (v,n+1)) \; \text{ and } \; 
\E_p(\cW_n) \ge \E_p(\cW_{n+1}).
$$  
\end{Thm}

\begin{Rem} Monotonicity of patterns. \label{Rem:patterns}
\begin{itemize}
\item 
For small graphs the onset of monotonicity occurs relatively early. 
Let $N(G)$ denote the minimal value of $N$ for which 
\eqref{equ:monoNp} holds. 
Then CAS-aided calculations that will be described in Subsection \ref{subsec:CAS}
give 
\begin{equation} \label{equ:table}
\begin{tabular}{c|c|c|c}
$G$ & $C_2$ & $C_3$ & $C_4$ \\ 
\hline 
$N(G)$ & $2$ & $2$ & $4$ 
\end{tabular}
\end{equation}
These calculations also show 
that Conjecture \ref{Con:layer} is true for the above graphs. 
\item 
It would be interesting to obtain explicit upper bounds on $N(G)$ for given graphs $G$. 
While we don't see how to extract this information from our proof of Theorem \ref{Thm:pattern} 
given in Section \ref{Sec:patterns}, 
we will show how to do that in a companion paper (\cite{KR})
making use of explicit bounds on the rate of convergence to the quasi-stationary distribution 
of our Markov chain (see \cite{CV}). 
With some effort we obtain $N(C_k) \le c 2^k$ for some constant $c$.
\end{itemize}
\end{Rem}

Concerning the monotonicity of the expected number of infected vertices per layer, 
we note that in view of Conjecture \ref{Con:layer} and  Proposition \ref{Prop:monodiff} 
we expect that 
\begin{displaymath} \label{equ:monoexp}
\forall n \ge 0 \forall p \in [0,1]: \E_p(\cW_{n}) \ge \E_p(\cW_{n+1}).  
\end{displaymath}
It is not too difficult to see that the above theorem implies that this holds for all $n \ge 0$ 
and sufficiently small $p$. In fact this can be seen more explicitly (and for graphs that are not necessarily finite) by a recursive approach, 
which we present as a complementary result: 
\begin{Thm} \label{Thm:number}
Let $G = (V,E)$ be a connected graph 
of bounded degree such that $\deg(x) \le \De$ for all $x \in V$.   
Then for percolation on $G^\times$ we have 
\begin{equation} \label{equ:numbermono}
\forall n \ge 0, p \in [0,\frac 1 {\De + 1.4}]: \E_p(\cW_n) \ge \E_p(\cW_{n+1}). 
\end{equation}
\end{Thm}
\begin{Rem} \label{Rem:number}
Monotonicity of the expected number of infected vertices. 
\begin{itemize}
\item 
The bound on $p$ in \eqref{equ:numbermono} can be improved
in particular cases. 
In case of $G = C_k$ and $G = \Z$ we have $\De = 2$, so 
the bound in \eqref{equ:numbermono} is $p \le \frac 1 {3.4} \approx 0.29$,  
but with some effort (and combinatorial arguments) 
one can show that $p \le 0.35$ suffices, see \cite{KR}. 
\item 
For further monotonicity results for the expectations we refer to \cite{Ri1}: 
If $G$ is locally finite we have $\E_p(\cW_n) = \infty$ for some $n$ iff 
$\E_p(\cW_n) = \infty$ for all $n$ (see Lemma 4), 
and under suitable symmetry assumptions on $G$ 
(that are satisfied in case of $G = C_k$ or $G = \Z$) 
we have $\E_p(\cW_0) \ge \E_p(\cW_n)$ for all $n \ge 0$ and $p \in [0,1]$ (see Corollary 2(a)).
\end{itemize}
\end{Rem}

\section{Comparing monotonicity properties}  \label{Sec:comparison}

\subsection{Comparing the monotonicity conjectures}

Here we prove Proposition \ref{Prop:bunkbed}. 
The main idea for (a) is that a suitable truncation 
of a layered graph is a bunkbed graph. 
For a detailed proof let $p \in (0,1)$.
We assume the validity of the bunkbed conjecture. 
Let $G$ be a finite connected graph, $o,v \in V(G)$ and $n \ge 0$. 
We note that for proving \eqref{equ:layer} it suffices to 
consider $n' = n+1$. 
For $m \ge n$ we let the subgraph of $G^\times$ induced by the vertex set 
$V(G) \times \{n-m,...,n+m+1\}$ be denoted by $G^{(m)} = (V^{(m)},E^{(m)})$, 
and we write $\lra_{(m)}$ for the connectivity relation induced 
by percolation on $E^{(m)}$, i.e. $u' \lra_{(m)} v'$ 
iff there is a path from $u'$ to $v'$ consisting of open bonds from $E^{(m)}$. 
We note that $G^{(m)}$ may be identified with 
the bunkbed graph $\tilde{G}^{(m)}_T$, 
where $\tilde{G}^{(m)}$ is defined to be the subgraph 
of $G^\times$ induced by the vertex set 
$G \times \{n-m,...,n\}$ and $T = G \times \{n\}$. 
Thus the validity of the bunkbed conjecture for $\tilde{G}^{(m)}_T$ implies  
$$
\P_p((o,0) \lra_{(m)} (v,n)) \ge 
\P_p((o,0) \lra_{(m)} (v,n+1)). 
$$
Since $\{(o,0) \lra_{(m)} (v,n)\} \uparrow 
\{(o,0) \lra (v,n)\}$ for $m \to \infty$ the above implies that 
$\P_p((o,0) \lra (v,n)) \ge 
\P_p((o,0) \lra (v,n+1))$ as desired. 

\bigskip 

The main idea for (b) is that $\Z^d$ can be considered to be 
a layered graph $(\Z^{d-1})^\times$ w.r.t. every direction, 
and $\Z^{d-1}$ can be approximated by $C_k^{d-1}$ with large $k$. 
For a detailed proof let $p \in (0,1)$ and $d \ge 2$. 
We assume the validity of the monotonicity conjecture  for the layered graphs 
$(C_k^{d-1})^\times$ for every $k \ge 2$. 
We note that for proving \eqref{equ:Zd} it suffices to 
show that for all $o \le w \in \Z^{d-1}$ and $n \ge 0$ 
we have $\P_p((o,0) \lra (w,n)) \ge \P_p((o,0) \lra (w,n+1))$, 
since we can go from any $u \ge o$ to any $v \ge u$ by successively increasing some coordinate by $1$. 
In order to infer this inequality from the given monotonicity property, it suffices 
to show that we have $\P_p(A_k) \to \P_p(A)$ for $k \to \infty$, 
where 
$$
A_k := \{o' \lra w' \text{ w.r.t. } (C_{2k+1}^{d-1})^\times\} \quad 
 \text{ and } \quad A := \{o' \lra w' \text{ w.r.t. } \Z^d\}  
$$
with $w' \in \Z^d$ and $o' = (o,0)$. 
(Here we label the vertices of $C_{2k+1}$ as $\{-k,...,k\}$
so that $o',w' \in C_{2k+1}^{d-1} \times \Z$ for sufficiently large $k$.) 
To show this convergence we consider the subgraph $G_{(k)} = (V_{(k)}, E_{(k)})$ of $\Z^d$ induced by the vertex set 
$\{-k,...,k\}^{d}$. 
$G_{(k)}$ naturally can also be considered a subgraph 
of $(C_{2k+1}^{d-1})^\times$. 
Let $\partial V_{(k)} = \{v' \in V_{(k)}: \exists i: |v_i'| = k\}$ 
be the boundary of $V_{(k)}$ and let $\lra_k$ be the connectivity
relation induced by percolation on $E_{(k)}$. 
We set
$$
B_k := \{\partial V_{(k)} \lra_k o' \not \lra_k w' \lra_k \partial V_{(k)}\}
$$
and note that $B_k^c \cap \{o' \lra w'\} =  \{o' \lra_k w'\}$. 
Thus both $B_k$ and $B_k^c \cap \{o' \lra w'\}$ only depend on edges in $E_{(k)}$, i.e. their probabilities w.r.t. percolation on $\Z^d$ or percolation 
on $(C^{d-1}_{2k+1})^\times$ are the same.  
This implies that $|\P_p(A_k) - \P_p(A)| \le $
\begin{align*}
\P_p(A_k \cap B_k) + |\P_p(A_k \cap B_k^c) - \P_p(A \cap B_k^c)| + \P_p(A \cap B_k)
\le 2 \P_p(B_k). 
\end{align*}
For the last term we note that (interpreting $B_k$ 
in terms of percolation on $\Z^d$) $B_k \downarrow B := \bigcap_k B_k$, 
and on $B$ the clusters of $o'$ and $w'$ 
are infinite and do not intersect, thus $\P_p(B) = 0$ 
by Burton-Keane (e.g. see [G]). 
So letting $k \to \infty$ we obtain $|\P_p(A_k) - \P_p(A)| \to 0$
as desired.  \qed

\subsection{Comparing monotonicity properties for layered graphs}

Here we prove Proposition \ref{Prop:monodiff}. 
Let $n \ge 0$ and $p \in (0,1)$. 
(i) $\Ra$ (ii) follows from the way connectivity events 
can be obtained from infection patterns.
For $n \ge 0$, $x \in M$ and $v \in V$ we write 
$* \lra_{x,n} v$ iff $(v,n)$ is infected from the pattern $x$ 
at layer $n$ via bonds in $E_{n+1..}$, i.e. 
iff there are $v_0',...,v_m' \in V \times \Z$ such 
that $v_0' = (u,n)$ for some $u \sim_x *$, $v_m' = (v,n)$ 
and for every $i$ either $v_i'v_{i+1}' \in E_{n+1..}$ such 
that $\cZ_{v_i'v_{i+1}'} = 1$ 
or $v_i' = (v_i,n), v_{i+1}' = (v_{i+1},n)$ for some $v_i \sim_x v_{i+1}$. 
For all $v \in V$ we have 
\begin{align*}
&\P_p((o,0) \lra (v,n)) = \sum_{x \in M^*} \P_p(\cX_n = x,(o,0) \lra (v,n)) \\
&= \sum_{x \in M^*} \P_p(\cX_n = x, * \lra_{x,n} v) 
=  \sum_{x \in M^*} \P_p(\cX_n = x)\P_p(* \lra_{x,n} v). 
\end{align*}
In the first step we have used that every open path from $(o,0)$ to $(v,n)$
contains an open subpath from $(o,0)$ to $V_n$ via bonds in $E_{..n}$. 
In the second step we have used that every path from $(o,0)$ to $(v,n)$ 
can be decomposed into a path from $(o,0)$ to $V_n$ via bonds in $E_{..n}$,
paths from $V_n$ to $V_n$ via bonds in  $E_{n+1..}$
and paths from $V_n$ to $V_n$ via bonds in $E_{..n}$, 
and that the connectivity of any pair of 
vertices in $V_n$ via bonds in $E_{..n}$
is encoded in $\cX_n$. 
In the third step we have used independence of the two events, 
which depend on disjoint sets of bonds. 
Repeating the above arguments we get a corresponding decomposition for 
$\P_p((o,0) \lra (v,n+1))$. 
Since $\P_p(* \lra_{x,n} v) = \P_p(* \lra_{x,n+1} v)$ by the 
translation invariance of percolation and 
$\P_p(\cX_n = x) \ge \P_p(\cX_{n+1} = x)$ for all $x \in M^*$ 
by (i) we thus see that $\P_p((o,0) \lra (v,n)) \ge \P_p((o,0) \lra (v,n+1))$.

\bigskip 

(ii) $\Ra$ (iii) is a straightforward consequence of the linearity of expectation: We have 
$\E_p(\cW_n) = \sum_{v \in V} \P_p((o,0) \lra (v,n))$
and similarly for $\E_p(\cW_{n+1})$. \qed

\section{Markov chain for infection patterns} \label{Sec:Markov}

\subsection{Markov property for infection patterns}

In all of this subsection we consider percolation on $G^\times$ with some parameter $p \in (0,1)$, 
where $G$ is a fixed finite connected graph. 
It will be useful to describe percolation events between two layers that start from a prescribed pattern:

\begin{Def} 
Let $0 \le k \le n$ and $y \in M$.
Let $\lra_{y,k,n}$ denote the connectivity relation induced by 
percolation on $E_{k+1..n}$
and the pattern $y$ in layer $V_k$. 
Let $\cX^{y,k}_{n}$ be the pattern in $V_n$ 
induced by $\lra_{y,k,n}$ with the infection originating in $y$, and let $\cX^y_n := \cX^{y,0}_n$. 
\end{Def} 
Writing out the above definition we have 
$u' \lra_{y,k,n} v'$ 
iff there are $v_0',...,v_m' \in V \times \Z$ such that 
$v_0' = u', v_m' = v'$ and for every $i$ either $v_i'v_{i+1}' \in E_{k+1..n} $ with $\cZ_{v_i'v_{i+1}'} = 1$ or $v_i' = (v_i,k), v'_{i+1} = (v_{i+1},k)$ with $v_i \sim_y v_{i+1}$. 
Furthermore we have $\cX^{y,k}_{n} = x$ iff   
\begin{align*}
&\forall u,v \in V: u \sim_x v \Lra (u,n) \lra_{y,k,n} (v,n) \quad \text{ and } \\ 
&\forall u \in V: u \sim_x * \Lra \exists w \in V: w \sim_y *, 
(u,n)  \lra_{y,k,n} (w,k). 
\end{align*}
In simple words, the event $\{\cX^{y,k}_n = x\}$ 
describes all configurations of open bonds in $E_{k+1..n}$ such that a prescribed pattern $y$ at layer $k$ produces the pattern $x$ at layer $n$. 
We also note that $\cX_n^{y,n} = y$ and in particular 
$\cX_0^y  = y$.  
\begin{Lem} 
For all $n \ge 0$ and $y,x_0,..,x_n \in M$ 
we have 
\begin{align} \label{equ:MarkovP}
&\P_p(\forall 0 \le i \le n: \cX_i = x_i) = 
\P_p(\cX_0 = x_0) \prod_{1 \le i \le n} \P_p(\cX_1^{x_{i-1}} = x_i) \text{ and }\\
\label{equ:MarkovPz}
&\P_p(\forall 0 \le i \le n: \cX_i^y = x_i) = 
\de_y(x_0) \prod_{1 \le i \le n} \P_p(\cX_1^{x_{i-1}} = x_i). 
\end{align} 
In particular $(\cX_n)_{n \ge 0}$ and 
$(\cX_n^y)_{n \ge 0}$ are Markov chains with the same transition probabilities $\pi_p(x,x') := \P_p(\cX_1^x = x')$ and initial distributions $\al_p := \P_p(\cX_0 = .)$ and $\de_y$ respectively. 
\end{Lem}

\Pf \eqref{equ:MarkovP} follows by noting that we have 
$$
\{\forall 0 \le i \le n: \cX_i = x_i\} = \{\cX_0 = x_0, 
\forall 1 \le i \le n: \cX_i^{x_{i-1},i-1} = x_i\}  
$$
by definition,
the random variables $\cX_0, \cX_1^{x_0,0},...,\cX_n^{x_{n-1},n-1}$ 
are independent since they depend on disjoint bond sets, 
and  $\cX_i^{x_{i-1},i-1} \sim \cX_1^{x_{i-1},0}$, 
which is a consequence of the translation invariance of percolation.  
\eqref{equ:MarkovPz} can be seen similarly.
\qed

\bigskip

In particular the above proves Proposition \ref{Prop:Markov}. 
We next collect some simple properties of the above Markov chain 
of patterns for later reference.
Let 
\begin{equation}
\xda := \{\{*\}\} \cup \{\{v\}: v \in V\} \quad \text{ and } \quad 
\xst := \{\{*\} \cup V\}
\end{equation}
denote the pattern without infection and connections and the pattern, where
everything is connected and infected, respectively.

\begin{Lem} \label{lem:transitions}
For the above Markov chain of patterns we have: 
\begin{align} 
\label{equ:positivity}
&\forall p,p' \in (0,1), x,y \in M: \pi_p(x,y) > 0 \Lra \pi_{p'}(x,y) > 0. \\
\label{equ:polynomial}
&\forall x,y \in M: \pi_p(x,y) \text{ is a polynomial function in $p$}. \\
\label{equ:absorbing}
&
\forall p \in (0,1), x \in \Mda, y \in M^*: \pi_p(x,y) = 0. \\
\label{equ:aperiodic}
&\forall p \in (0,1),x \in M: \pi_p(x,x) > 0. \\
\label{equ:nobonds}
&\forall p \in (0,1),x \in M: \pi_p(x,\xda) > 0. 
\end{align}
\end{Lem}

\Pf For \eqref{equ:positivity} we note that a specific configuration of bonds 
in a certain layer has positive probability w.r.t. $\P_p$ if it has 
positive probability w.r.t. $\P_{p'}$. 
For \eqref{equ:polynomial} we similarly note that a specific configuration of bonds 
in a certain layer has a probability w.r.t. $\P_p$ of the form $p^k(1-p)^l$.  
\eqref{equ:absorbing} is immediate from the definition of the Markov chain. 
\eqref{equ:aperiodic} can be seen by considering a layer, where all 
vertical bonds are open and all horizontal bonds are closed. 
Similarly, \eqref{equ:nobonds} can be seen by considering a layer, where all 
bonds are closed. 
\qed

\bigskip 

As usual, for $x,y \in M$ we write $x \to y$ if $y$ can be reached from $x$, i.e. for some $n$ we have $\pi_p^n(x,y) > 0$. 
We note that by \eqref{equ:positivity} this does not depend on the value 
of $p \in (0,1)$. 
By \eqref{equ:absorbing}, $\Mda$ is an absorbing set of states.
The decomposition of $M$ into communicating classes can be quite 
complicated and in particular for most graphs the set $\Mda$ will not be irreducible. 
However, if it is irreducible, then \eqref{equ:aperiodic} implies aperiodicity. 
Finally, in preparation for the proof of Theorem \ref{Thm:pattern} we note 
the following consequences of the Markovian structure. 

\begin{Lem} \label{lem:monoN} For the above Markov chain of patterns we have:  \nobr 
\begin{itemize}
\item[(a)] Let $n \ge 0$. If for all $y,x \in M^*$ we have 
 $\P_p(\cX_n^y = x) \ge \P_p(\cX_{n+1}^y = x)$, then 
 for all $x \in M^*$  $\P_p(\cX_n = x) \ge \P_p(\cX_{n+1} = x)$. 
\item[(b)] Let $N \ge 0$. If for all $y,x \in M^*$ we have 
$\P_p(\cX_N^y = x) \ge \P_p(\cX_{N+1} ^y= x)$, 
then for all $n \ge N$ and all $y,x \in M^*$  
$\P_p(\cX_n^y = x) \ge \P_p(\cX_{n+1}^y = x)$. 
\end{itemize}
\end{Lem}

\Pf For (a) it suffices to note that for $x \in M^*$
$$
\P_p(\cX_n = x) = \sum_{y \in M^*} \P_p(\cX_0 = y)\P_p(\cX_n^y = x)
$$
and similarly for $n+1$ instead of $n$. Here we have used that 
$\P_p(\cX_n^y = x) = 0$ for $y \in \Mda, x \in M^*$ by \eqref{equ:absorbing}. 
For (b) it suffices to note that for all $x,y \in M^*$ 
\begin{align*}
\P_p(\cX_{n}^y = x) &= \sum_{x' \in M^*} \P_p(\cX_{N}^y = x') 
\P_p(\cX_{n-N}^{x'} = x) 
\end{align*}
and similarly for $n+1$ and $N+1$ instead of $n$ and $N$, 
where again we have used \eqref{equ:absorbing}.  \qed 

%

\subsection{CAS computations for small graphs} \label{subsec:CAS}

In this subsection we will outline how to obtain Table \eqref{equ:table}, i.e. 
on how to obtain optimal bounds on the onset of monotonicity with CAS-assistance. 
Similar calculations are possible for any sufficiently small graph $G$. 
In some cases, e.g. for $G = C_2 = L_2$ (for which $G^\times$ is the bi-infinite ladder graph), 
these computations can be carried out by hand. 
We will use this particular example to illustrate the general procedure. 

Since not all partitions of $V$ can be connectivity patterns of percolation on $G^\times$, 
let us first reduce the state space of the Markov chain.
We consider a version $\cX_n^-$ of the Markov chain that only records the connectivity in the layer $V_n$, 
but not which vertices are infected, i.e. $\cX_n^- \in M^\da$ 
and $\cX_n^- = x$ for some $x \in M^\da$ iff 
$\forall u,v \in V: u \sim_x v \Lra (u,n) \lra_{\le n} (v,n)$.  
From the proof of Proposition \ref{Prop:Markov} it can be seen, 
that $\cX_n^-$ indeed is a Markov chain with transition probabilities $\pi_p^-(x,y) = \pi_p(x,y)$ for all $x,y \in M^\da$. $\cX_n^-$ is stationary, which follows from its definition and the translation invariance of percolation on $G^\times$. 
Since a.s. for some $k \le 0$ all bonds of $E_k$ are closed, 
$\cX_n^-$ only takes values in $M^- := \{x \in M^\da: x_\da \to x\}$. Since 
$x \to x_\da$ for all $x \in M^\da$ by \eqref{equ:nobonds}, 
the Markov chain $\cX_n^-$ is irreducible.  
In case of $G = C_2$ we have 
$M^- = \{\{1\},\{2\},\{*\}\},\{\{1,2\},\{*\}\}$ and 
$$
\pi^-_p = \begin{pmatrix}
1-p & p \\ 
(1-p)(1-p^2) & p(1+p-p^2)
\end{pmatrix}. 
$$
By definition we have $\cX_n^- = f_\da (\cX_n)$, where 
$f_\da: M \to M^\da$ deletes the infection from a given partition, 
i.e. $f_\da(x) := \{\{*\}\} \cup \{A \stm \{*\}\!: A \in x\}$. 
In particular $\cX_n$ only takes values in $f_\da^{-1}(M^-)$. 
Since we are only interested in states with infections, we can further reduce the state space by combining all uninfected states into a single state $\da$. 
We thus consider $\cX_n' := f_*(\cX_n)$, 
where  $f_*: M \to M^* \cup \{\da\}$ is the identification map 
$f_*(x) := x$ for $x \in M^*, f_*(x) := \da$ for $x \in M^\da$. 
$\cX_n'$ is a Markov chain with state space 
$M' = (f_\da^{-1}(M^-) \cap M^*) \cup \{\da\}$ and transition 
probabilities $\pi'_p(x,y) = \pi_p(x,y)$ for all $x,y \in M' \stm \{\da\}$. 
In case of $G = C_2$ we have  
$M' = \{\da, \{\{1,*\},\{2\}\}, \{\{2,*\},\{1\}\}, \{\{1,2,*\}\}\}$ and 
$$
\pi'_p = \begin{pmatrix}
1 & 0 & 0 & 0 \\ 
1-p & p(1-p) & 0 & p^2\\
1-p & 0 & p(1-p) & p^2 \\ 
(1-p)^2 & p(1-p)^2 & p(1-p)^2 & p^2(3-2p)
\end{pmatrix}.
$$
In general, the transition probabilities of $\cX_n^-$ and $\cX_n'$ will be polynomials in $p$ of degree $ \le |E(G)| + |V(G)|$. 
The initial distribution $\al'_p(x) = \P_p(\cX_0' = x)$ of the chain $\cX_n'$ can be obtained from the chain $\cX_n^-$: 
Since $\cX_n^-$ is irreducible and stationary, 
the distribution of $\cX_0^-$ is the unique stationary 
distribution $\al_p^-$ of this chain. 
It can be computed as a left eigenvector of $\pi_p^-$ 
corresponding to the eigenvalue $1$, e.g. using Gaussian elimination. 
Since $\cX_0^- = f_\da(\cX_0')$ and $o \sim_{\cX_0'} *$,  
we then can set 
$\al'_p(x) := 0$ if $o \not \sim_x *$ and 
$\al'_p(x) := \al_p^-(f_\da(x))$ if $o \sim_x *$. 
In case of $G = C_2$ we obtain 
$$
\al_p^- = \frac 1 {c_p}(((1-p)(1-p^2),p)
\quad \text{ and } \quad \al'_p = \frac 1 {c_p}(0,((1-p)(1-p^2),0,p), 
$$ 
where $c_p$ is a normalizing constant. 
By Gaussian elimination in general the entries of $\al_p^-$ are rational functions of 
$p$, and by irreducibility all entries are strictly positive 
for all $p \in (0,1)$. 
Thus we can always choose $c_p$ such that $c_p >0$ for all $p \in (0,1)$
and the entries of $c_p \al_p^-$ (and thus also of $c_p \al_p'$) 
are polynomials in $p$. 
In view of (b) of \autoref{lem:monoN} it suffices 
to find an $N$ such that  all entries of $A(p,N) := c_p \al_p' {\pi_p'}^N - c_p \al_p' {\pi_p'}^{N+1}$ are positive for all $p \in (0,1)$. 
The positivity of these polynomials on $(0,1)$ can be decided using Sturm's method (e.g. see \cite{BPR}).
We thus are left to consider $A(p,n)$ for $n \in \{0,1,2,...\}$ 
until for $n = N$ $A(p,n)$ has the desired property. 

\bigskip 

If $N$ is not too large, we can thus also check that 
$\P_p((o,0) \lra (v,n)) \ge \P_p((o,0) \lra (v,n+1))$ for all 
$v \in V$, $p \in [0,1]$ and $n \ge 0$. 
For $n \ge N$ this follows directly from Proposition \ref{Prop:monodiff}, and 
the cases $n \in \{0,...,N-1\}$ can be checked separately using a similar idea to the one used in the proof of Proposition~\ref{Prop:monodiff}: 
Let $\cX_n$ denote (as before) the pattern in $V_n$ induced by percolation on $E_{..n}$. 
Let $\cY_{n+1}$ denote the uninfected pattern in $V_{n+1}$ 
induced by percolation on $E_{n+1}^h \cup E_{n+2..}$.  
Let $\cZ_{n+1}^v$ denote the configuration of bonds in $E_{n+1}^v$. 
We note that $\cX_n, \cZ_{n+1}^v, \cY_{n+1}$ are independent, 
$\cX_n \sim \al'_p{\pi_p'}^n$, $\cZ_{n+1}^v$ is a Bernoulli sequence with parameter $p$, and $\cY_{n+1} \sim \al_p^-$ by reflection invariance of percolation on $G^\times$. 
Let $A(v,n)$ denote the set of all triples $(x,y,z)$, such that 
$x$ is an infected pattern, $y$ is an uninfected pattern, $z$ is a bond configuration of a layer of vertical bonds,  
and there is a path $v_0',...,v_m' \in V_n \cup V_{n+1}$ 
such that $v_0' = (v_0,n)$ with $v_0 \sim_x *$, $v_m' = (v,n)$ and for each $i$ either $v_i'= (v_i,n), v_{i+1}' = (v_{i+1},n)$ 
with $v_i \sim_x v_{i+1}$, 
or $v_i'v_{i+1}' \in E_{n+1}^v$ with $z_{v_i'v_{i+1}'} = 1$ 
or $v_i'= (v_i,n+1)$, $v_{i+1}' = (v_{i+1},n+1)$ 
with $v_i \sim_y v_{i+1}$. Then we have 
\begin{align*}
\P_p((o,0)\lra (v,n)) &= \sum_{(x,y,z) \in A(v,n)} 
\P_p(\cX_n = x, \cY_{n+1} = y, \cZ_{n+1}^v = z)\\ 
&= \sum_{(x,y,z) \in A(v,n)} 
\al_p'{\pi_p'}^n(x) \al_p^-(y) \P_p(\cZ_{n+1}^v=z)
\end{align*}
and a similar expression for $(v,n+1)$. 
The probabilities $\al_p'{\pi_p'}^n(x)$ and  $\al_p^-(y)$ 
already have been computed while determining $N$, 
$\P_p(\cZ_{n+1}^v = z) = p^{k_z}(1-p)^{l_z}$ 
can easily be computed, and so 
$c_p^2\P_p((o,0)\lra (v,n))-c_p^2\P_p((o,0)\lra (v,n+1))$
is a polynomial in $p$ that can readily be computed and its positivity 
for $p \in (0,1)$ can be checked using Sturm's method as above.

\begin{Rem} CAS computations for small graphs $G$. 
\begin{itemize}
\item 
We have implemented the above ideas (in Python) to compute $N$
and to verify \eqref{equ:layer} in case of $G = C_k$ 
for $k \in \{2,3,4\}$. 
\item 
While it may be possible to extend these calculations to cover $k = 5$, for larger values of $k$ theses CAS computations no longer are feasible since the state spaces for the Markov chains are too big. We note that for $k = 5$ we have $|M'| = 127$ and $|M^-| = 42$ and the entries of $\pi_p'$ and $\pi_p^-$ are polynomials of degree 10. 
\end{itemize}
\end{Rem}

\section{Partition patterns per layer} \label{Sec:patterns}

Here we prove Theorem \ref{Thm:pattern}. 
We proceed by considering small, intermediate and large 
values of $p$ separately.

\subsection{Intermediate values of $p$}

The main ingredient here is a general monotonicity result 
for Markov chains which relies on 
asymptotics for transition probabilities. 
As reference for the needed asymptotics we use \cite{Li1}, 
for more general results see \cite{Li2} and the references therein.  
While the result of \cite{Li1} is purely algebraic, 
we would like to point out that it allows for a probabilistic interpretation in terms of quasi-stationary distributions, see 
\cite{DP}.  

\begin{Prop} \label{prop:markovmono}
Let $\pi$ be the transition matrix of an arbitrary Markov chain 
on a finite state space $M$. 
Suppose the chain is locally aperiodic, i.e. for every $x \in M$ 
either $\pi^n(x,x) = 0$ for all $n \ge 1$ or there is an $N$ such that 
$\pi^n(x,x) > 0$ for all $n \ge N$. Let $y,x \in M$. 
\begin{itemize}
\item[(a)]
Either there is an $N$ such that $\forall n \ge N: \pi^n(y,x) = 0$ or there are constants 
$\la > 0$, $\si \in \{1,2,...\}$ and $\ga > 0$ such that 
\begin{equation} \label{equ:markovasy}
\pi^n(y,x) \sim \ga n^{\si - 1} \la^n \quad \text{ for } n \to \infty. 
\end{equation}
\item[(b)] 
Suppose that $x$ is transient.
Then there is an $N$ such 
\begin{equation} \label{equ:markovmono}
\forall n \ge N: \pi^n(y,x) = 0 \quad \text{ or } \quad 
\forall n \ge N: \pi^{n+1}(y,x) < \pi^n(y,x). 
\end{equation}
\end{itemize}
\end{Prop}

\Pf For a proof of (a)  we refer to \cite{Li1} (Theorem, page~6). 
We note that this proof is purely algebraic 
and gives fairly explicit descriptions of the constants $\la,\si,\ga$.
For a proof of (b) let $x$ be transient. 
If $\forall n \ge N: \pi^n(y,x) = 0$ we are done. Thus by (a) we may assume that \eqref{equ:markovasy} holds, which implies  
$$
\frac{\pi^{n+1}(y,x)}{\pi^n(y,x)} \to \la \quad \text{ for } n \to \infty. 
$$
It suffices to show that $\la < 1$, which follows readily from the transience of $x$:  
The expected number of visits in $x$ 
from a chain started in $y$ is finite, i.e. $\sum_n \pi^n(y,x) < \infty$, 
which implies $\pi^n(y,x) \to 0$ for $n \to \infty$. 
By $\eqref{equ:markovasy}$ this implies that $\ga n^{\si - 1} \la ^n \to 0$, 
and thus $\la < 1$.
 \qed 
 
\bigskip 
 
Now we return to the setting of Theorem \ref{Thm:pattern} with a fixed graph $G$.

\begin{Prop} \label{prop:intermediate}
For all $0 < p_0 < p_1 < 1$ there is an $N \ge 0$ such that 
\begin{equation}\label{equ:intermediate}
\forall n \ge N \forall p \in [p_0,p_1] \forall y,x \in M^*: 
\P_p(\cX_n^y =x) \ge \P_p(\cX_{n+1}^y = x). 
\end{equation}
\end{Prop}

\Pf We first consider a fixed $p' \in [p_0,p_1]$. 
The state space $M$ of the considered Markov chain is finite.
The Markov chain is locally aperiodic since $\pi_{p'}^n(x,x) > 0$ for all 
$n \ge 0$ by \eqref{equ:aperiodic}. 
For every $x \in M^*$ we have $\pi_{p'}(x,\xda) > 0$ and 
$\pi_{p'}^n(\xda,x) = 0$ for all $n \ge 0$ by \eqref{equ:nobonds} and 
\eqref{equ:absorbing}, so $x$ is transient. 
Thus, by Proposition \ref{prop:markovmono} for all $y,x \in M^*$ there is an 
$N_{p',y,x}$ 
such that \eqref{equ:markovmono} holds for $\pi = \pi_{p'}$ and $N = N_{p',y,x}$. 
Since $M^*$ is finite, we can define $N_{p'} := \max\{N_{p',y,x}: y,x \in M^*\}$. By definition of $N_{p'}$ we have 
$$
\forall y,x \in M^*: \quad \pi_{p'}^{N_{p'}+1}(y,x) = \pi_{p'}^{N_{p'}}(y,x) = 0
\quad \text{ or } \quad \pi_{p'}^{N_{p'}+1}(y,x) < \pi_{p'}^{N_{p'}}(y,x). 
$$
By \eqref{equ:positivity} $\pi_{p'}^n(y,x) = 0$ implies that  
$\forall p \in (0,1): \pi_{p}^n(y,x) = 0$, 
and by \eqref{equ:polynomial} $\pi_{p}^n(y,x)$ is a continuous function of $p$. 
Thus the above implies that there is an open neighborhood $U_{p'}$
of $p'$ in $(0,1)$ such that 
$$
\forall p \in U_{p'} \forall y,x \in M^*: \pi_{p}^{N_{p'}+1}(y,x) \le \pi_{p}^{N_{p'}}(y,x). 
$$
By Lemma \ref{lem:monoN} this implies that 
$$
\forall p \in U_{p'} \forall y,x \in M^* \forall n \ge N_{p'}: 
\P_p(\cX_{n+1}^y = x) \le \P_p(\cX_{n}^y = x).  
$$
By compactness, $[p_0,p_1] \subset \bigcup_{i \in I} U_{p_i}$, 
for some finite set $I$ and some $p_i \in [p_0,p_1]$.  
Setting  $N_{p_0,p_1} := \max\{N_{p_i}: i \in I\}$ we have 
$$
\forall i \in I \forall p \in U_{p_i} \forall y,x \in M^* \forall n \ge N_{p_0,p_1}: 
\P_p(\cX_{n+1}^y = x) \le \P_p(\cX_{n}^y = x).  
$$
Thus \eqref{equ:intermediate} holds for $N = N_{p_0,p_1}$. \qed

\subsection{Small values of $p$} \label{Sec:Small}

In this case we use that going from a pattern $y$ to a pattern $x$ in a fixed large number of steps, 
in most layers  the number of open bonds will be as small as possible. 
With high probability there will be two consecutive layers with the same configuration of bonds, where all horizontal bonds are closed. 
By discarding one of these two layers we will be able to compare 
$\pi_p^{n+1}(y,x)$ and $\pi_p^{n}(y,x)$. 
To make this into a rigorous argument we introduce some notation and collect some observations. 

\begin{Def} Let  $\cN_n := |\{e \in E_n\!: \cZ_e \!=\! 1\}|$, for $y,x \in M^*$ such that $y \to x$ let 
\begin{align*}
&m_{y,x} := \min\{\cN_i(\om): \exists n \ge 1: 
1 \le i \le n, \om \in \{\cX_n^{y} = x\}\} \text{ and } \\
&l_{y,x} := \min\{n \ge 1: \exists 1 \le i \le n: 
\{\cX_n^{y} = x,\cN_i = m_{y,x}\} \neq \emptyset\}, 
\end{align*}
and let 
$l := \max\{l_{y,x}: y,x \in M^* \text{ such that } y \to x\}.$ 
\end{Def}
Thus $\cN_n$ is the number of open bonds in layer $E_n$, 
$m_{y,x}$ is the minimal number of open bonds in some layer of a bond configuration 
going from pattern $y$ to pattern $x$ in an arbitrary number of steps, 
$l_{y,x}$ is the minimal number of steps in which you can go from $y$ to $x$ and observe a layer with the minimal number of open bonds in between. 
We note that $m_{y,x} \ge 1$, $l_{y,x} < \infty$, $l < \infty$ 
and for every $n \ge l$ and all $y,x \in M^*$ such that 
$y \to x$ it is possible to go from $y$ to $x$ in 
$n$ steps and observe a layer with the minimal number of open bonds in between (making use of \eqref{equ:aperiodic}).
Since one might speculate whether values $m_{y,x} > 1$ are possible at all, let us provide a simple example: 
For $G = L_3$ and $x = y = \{\{*,1\},\{0,2\}\}$ 
we have $m_{y,x} = 3$.

\begin{Lem} \label{lem:minopen} Layers with the minimal number of open bonds.
For all $y,x \in M^*$ such that $y \to x$ we have 
\begin{align}
\label{equ:minopen1}
&\forall 1 \le i \le n: \{\cX_n^{y}= x, \cN_i = m_{y,x}\} \subset \{\forall e \in E_{i}^h: \cZ_e = 0\}, \\
\label{equ:minopen2}
&\forall 1 \le i < n: \{\cX_n^{y}= x, \cN_i = \cN_{i+1} = m_{y,x}\} \subset 
\{\forall e \in E_i: \cZ_e = \cZ_{\tau(e)}\}, \end{align}
i.e. 
in a layer with the minimal number of open bonds  all horizontal bonds are closed, and 
in two consecutive layers with the minimal number open bonds, 
the open bonds match up.  
\end{Lem}

\Pf For \eqref{equ:minopen1} let $z = (z_e)_{e \in E_{1..n}}$ be a bond configuration 
contributing to the event $\{\cX_n^{y}= x, \cN_i = m_{y,x}\}$. 
Let $\bz = (\bz_e)_{e \in E_{1..n+1}}$ be the bond configuration 
obtained from $z$ by extending the open vertical bonds 
in $E_i$, i.e. $\bz_e = z_e$ for $e \in E_{1..i-1} \cup E_{i}^v$, 
$\bz_e = 0$ for $e \in E_{i}^h$, 
$\bz_e = z_{\tau^{-1}(e)}$ for $e \in E_{i+1..n+1}$. 
By construction $\bz$ contributes to 
$\{\cX_{n+1}^{y} = x,\cN_i = m, \cN_{i+1} = m_{y,x}\}$ for some $m$. 
If $z$ has an open bond in $E_{i}^h$, i.e. $\bz$ has an open bond in $E_{i+1}^h$, then $m < m_{y,x}$ 
by construction, contradicting the minimality of $m_{y,x}$. 
\\
For \eqref{equ:minopen2} let $(z_e)_{e \in E_{1..n}}$ be a bond configuration 
contributing to the event $\{\cX_n^{y}= x, \cN_i = \cN_{i+1} = m_{y,x}\}$. 
Suppose that the open bonds in the two layers do not match up. 
Since the number of open bonds in these layers is the same 
and the horizontal bonds match by \eqref{equ:minopen1}, 
there must be an open bond in $E_{i+1}^v$ 
such that the corresponding bond in $E_i^v$ is closed. 
By closing this open bond we obtain a bond configuration 
$(z'_e)_{e \in E_{1..n}}$, which contributes to 
$\{\cX_n^{y}= x, \cN_{i+1}= m_{y,x}-1\}$, 
contradicting the minimality of $m_{y,x}$. 
\qed

\begin{Prop} \label{prop:smallprep}
There are constants $c_1,c_2,c_3 \ge 1$ depending only 
on $G$ such that for all $n \ge l$, $p \in (0,\frac 1 2)$ and 
$y,x \in M^*$ such that $y \to x$ we have 
\begin{align} 
\label{equ:smallpmatching}
&\P_p(\cX_{n+1}^y = x, A_{y,x,n+1}) \le pn \P_p(\cX_n^y = x),\\
\label{equ:smallpnotmatching}
&\P_p(\cX_{n+1}^y = x, A_{y,x,n+1}^c) \le c_1^n p^{\frac n 2 (2m_{y,x}+1)},
\\
\label{equ:smallpprob}
&\P_p(\cX_n^y = x) \ge p^{c_2} c_3^{-n} p^{n m_{y,x}} ,  
\end{align}
where  $A_{y,x,n+1} := \{\exists 1 \le i \le \frac {n+1} 2: \cN_{2i-1} = m_{y,x} = \cN_{2i}\}$.
\end{Prop}

\Pf Let $b := |E_1|$. For \eqref{equ:smallpmatching} we note that 
by definition of $A_{y,x,n+1}$
$$
\P_p(\cX_{n+1}^y = x, A_{y,x,n+1})
\le \sum_{1 \le i \le \frac{n+1}2} \P_p(\cX_{n+1}^y = x, 
\cN_{2i-1} = m_{y,x} = \cN_{2i})
$$
using a union bound, 
and we note that for every $i$ 
$$
\P_p(\cX_{n+1}^y = x, 
\cN_{2i-1} = m_{y,x} = \cN_{2i})
\le \P_p(\cX_{n}^y = x) p.
$$
To see this, we consider a bond configuration  
$z = (z_e)_{e \in E_{1..n+1}}$ contributing to $\{\cX_{n+1}^{y} = x, 
\cN_{2i-1} = m_{y,x} = \cN_{2i}\}$. 
We decompose this configuration 
into $(z_e)_{e \in E_{2i-1}}$ and $\bz = (\bz_e)_{e \in E_{1..n}}$, 
which is obtained from $z$ be deleting layer $E_{2i-1}$, i.e. $\bz_e = z_e$ for $i \in E_{1..2i-2}$ 
and $\bz_e = z_{\tau(e)}$ for $e \in E_{2i-1..n}$.  
By the previous lemma $\bz$ contributes to $\{\cX_n^{y} = x\}$, 
and since $m_{y,x} \ge 1$ we have 
$$
\P_p(\forall e \in E_{2i-1}: \cZ_e = z_e) 
= p^{m_{y,x}} (1-p)^{b - m_{y,x}} \le p.
$$ 
This establishes the above estimate and using 
$\frac{n+1} 2 \le n$ we thus get \eqref{equ:smallpmatching}. 
For \eqref{equ:smallpnotmatching} we note that 
$$
\P_p(\cX_{n+1}^y = x, A_{y,x,n+1}^c)
\le \P_p(\forall 1 \le i \le \frac{n+1} 2: \cN_{2i-1} + \cN_{2i} > 2 m_{y,x})  
$$
by definition of $A_{y,x,n+1}$ and $m_{y,x}$. 
Furthermore for every $i > 0$ 
\begin{align*}
&\P_p(\cN_{2i-1} + \cN_{2i} > 2 m_{y,x}) 
= \sum_{2 m_{y,x} + 1 \le k \le 2b} \binom {2b} k p^k(1-p)^{2b-k}\\
&\le \sum_{0 \le k \le 2b} \binom {2b} k p^{ 2 m_{y,x} + 1}
= 2^{2b}p^{ 2 m_{y,x} + 1}.
\end{align*}
Combining this with the independence of $\cN_{2i-1} + \cN_{2i}$ $(i > 0)$, setting $c_1 := 2^{2b}$ and using $\frac{n+1} 2 \le n, \frac{n+1}2 \ge \frac n 2$
we obtain \eqref{equ:smallpnotmatching}. 
For \eqref{equ:smallpprob} we note that by definition of $l$,  
for some $i \in \{1,...,l\}$ we have 
$\{\cX_l^y = x, \cN_i = m_{y,x}\} \neq \emptyset$. 
Let  $z = (z_e)_{e \in E_{1..l}}$ be a bond configuration 
contributing to this event, 
and let $\bz = (\bz_e)_{e \in E_{1..n}}$ 
be the bond configuration obtained from $z$ 
by inserting $n-l$ additional layers above layer $E_i$, 
such that in these additional layers 
the open bonds precisely match the positions of the open bonds in $E_i$. 
By the above lemma $\bz$ contributes to $\{\cX_n^y =x\}$
and thus 
\begin{align*}
&\P_p(\cX_n^y = x) \ge \P_p(\forall e \in E_{1..n}: \cZ_e = \bz_e)\\ 
&= \P_p(\forall e \in E_{1..l}: \cZ_e = z_e) (p^{m_{y,x}}(1-p)^{b - m_{y,x}})^{n-l}. 
\end{align*}
Noting that $p \le \frac 1 2$ 
we estimate 
$\P_p(\forall e \in E_{1..l}: \cZ_e = z_e) \ge p^{bl}$, 
$(1-p)^{b-{m_{y,x}}} \ge (\frac 1 2)^b$
and $n-l \le n$. 
Thus we obtain \eqref{equ:smallpprob} by setting $c_2 := bl$ and $c_3 := 2^b$.
\qed 

\bigskip 

The above ideas enable us to treat the case
of small values of $p$: 

\begin{Prop} \label{prop:small}
There is a $p_0 \in (0,\frac 1 2)$ and an $N \ge 0$ such 
that 
\begin{equation}\label{equ:small}
\forall n \ge N \forall p \in (0,p_0] \forall y,x \in M^*: 
\P_p(\cX_n^y =x) \ge \P_p(\cX_{n+1}^y = x). 
\end{equation}
\end{Prop}

\Pf Choosing $c_1,c_2,c_3  \ge 1$ according to the above proposition,
we fix $N \ge \max\{l,4c_2\}$ and consider $p \le \frac 1 2$ and $y,x \in M^*$. In case of $y \to x$ we get 
\begin{align*}
&\P_p(\cX_{N+1}^y = x) = \P_p(\cX_{N+1}^y = x, A_{y,x,N+1}) + \P_p(\cX_{N+1}^y = x, A_{y,x,N+1}^c)\\
&\le \P_p(\cX_N^y = x) \big( p N + c_1^N p^{\frac N 2} c_3^{N} p^{-c_2}\big)
\le \P_p(\cX_N^y = x)
\big( p N + (c_1c_3)^N p^{\frac N 4}\big). 
\end{align*}
If $p$ is sufficiently small, we thus have 
$\P_p(\cX_{N+1}^y = x) \le \P_p(\cX_N^y = x)$. 
In case of $y \not \to x$ we have 
$\P_p(\cX_{N+1}^y = x) = 0 =  \P_p(\cX_N^y = x)$. 
In conclusion there is a value $p_0 \in (0,\frac 1 2)$ 
(depending on $G$ only) such that for all $p \in (0,p_0]$ and 
$y,x \in M^*$ we have 
$\P_p(\cX_N^y = x) \ge \P_p(\cX_{N+1}^y = x)$, and 
\eqref{equ:small} follows using Lemma \ref{lem:monoN}. 
\qed 

\subsection{Large values of $p$}\label{Sec:large}

In this case we will use that going from a pattern $y$ to a pattern $x$ in a fixed large number of steps, in most layers  the number of open bonds will be as large as possible. If it is not possible that all bonds in a layer are open, 
we can mimic the argument of the previous subsection: 
With high probability there are two consecutive layers 
with the same configuration of bonds, 
where all vertical bonds are open, and we proceed 
by discarding one of these layers. 

\begin{Def} Let $\cN_n' := |\{e \in E_n: \cZ_e= 0\}|$, for $y,x \in M^*$ such that $y \to x$ let 
\begin{align*}
&m'_{y,x} := \min\{\cN_i'(\om): \exists n \ge 1: 
1 \le i \le n, \om \in \{\cX_n^{y} = x\}\} \text{ and } \\
&l_{y,x}' := \min\{n \ge 1: \exists 1 \le i \le n: 
\{\cX_n^{y} = x,\cN_i' = m_{y,x}'\} \neq \emptyset\}, 
\end{align*}
and let $
l' := \max\{l'_{y,x}: y,x \in M^* \text{ such that } y \to x\}.
$ 
\end{Def}
Thus $\cN'_n$ is the number of closed bonds in layer $E_n$, 
$m_{y,x}'$ is the minimal number of closed bonds 
in some layer of a bond configuration going from pattern $y$ to pattern $x$ in an arbitrary number of steps, 
$l_{y,x}'$ is the minimal number of steps in which you can go from $y$ to $x$ and observe a layer with the minimal number of closed bonds in between. 
We note that $l_{y,x}' < \infty$, $l' < \infty$ 
and for every $n \ge l'$ and all $y,x \in M^*$ such that 
$y \to x$ it is possible to go from $y$ to $x$ in 
$n$ steps and observe a layer with the minimal number of closed bonds in between (making use of \eqref{equ:aperiodic}).
We stress that here $m_{y,x}' = 0$ cannot be ruled out.

\begin{Lem} \label{lem:minclosed} Layers with the minimal number of closed bonds.
For all $y,x \in M^*$ such that $y \to x$ we have 
\begin{align}
\label{equ:minclosed1}
&\forall 1 \le i \le n: \{\cX_n^y= x, \cN_i' = m_{y,x}'\} \subset \{\forall e \in E_i^v: \cZ_e = 1\}, 
\\
\label{equ:minclosed2}
&\forall 1 \le i < n: 
\{\cX_n^{y}= x, \cN'_i = \cN'_{i+1} = m'_{y,x}\} \subset 
\{\forall e \in E_i: \cZ_e = \cZ_{\tau(e)}\}, 
\end{align}  
i.e. 
in a layer with the minimal number of closed bonds  all vertical bonds are open, and 
in two consecutive layers with the minimal number closed bonds, 
the closed bonds match up.  
\end{Lem}

\Pf This is very similar to the proof of Lemma \ref{lem:minopen}. 
For \eqref{equ:minclosed1} let $z = (z_e)_{e \in E_{1..n}}$ 
be a bond configuration contributing to the event  
$\{\cX_n^y= x, \cN_i' = m_{y,x}'\}$. 
Let $\bz = (\bz_e)_{e \in E_{1..n+1}}$ be the bond configuration 
obtained from $z$ by inserting an additional layer 
in between layers $E_{i}$ and $E_{i+1}$ 
with all vertical bonds open and the horizontal bonds as in $E_i$, i.e. $\bz_e = z_e$ for all $e \in E_{1..i}$, $\bz_e = 1$ 
for all $e \in E_{i+1}^v$ and $\bz_e = z_{\tau^{-1}(e)}$ for 
all $e \in E_{i+2..n+1} \cup  E_{i+1}^h$. 
By construction $\bz$ contributes to 
$\{\cX_{n+1}^{y} = x,\cN_i' = m_{y,x}', \cN_{i+1}' = m'\}$ for some $m'$. 
If $z$ and thus $\bz$ has  a closed bond in $E_i^v$, 
then $m' < m_{y,x}'$ by construction, contradicting the minimality of $m_{y,x}'$. 
\\
For \eqref{equ:minclosed2} let $(z_e)_{e \in E_{1..n}}$ 
be a bond configuration contributing to the event
$\{\cX_n^{y}= x, \cN'_i = \cN'_{i+1} = m'_{y,x}\}$. 
Suppose that the closed bonds in the two layers do not match up. 
Since the number of closed bonds in these layers is the same 
and the vertical bonds match by \eqref{equ:minclosed1}, 
there must be a closed bond in $E_{i+1}^h$ such that the corresponding bond in $E_i^h$ is open. 
By opening this closed bond we obtain a bond configuration $(z_e')_{e \in E_{1..n}}$, which contributes to 
 $\{\cX_n^{y}= x, \cN_{i}'= m_{y,x}'-1\}$,
 and this contradicts the minimality of $m'_{y,x}$. 
\qed

\begin{Prop} \label{prop:largeprep}
There are constants $c_1,c_2,c_3 \!\ge \! 1$ depending only 
on $G$ such that for all $n \ge l'\!$, $p \in (\frac 1 2,1)$ and 
$y,x \in M^*$ such that $y \to x$ and $m_{y,x}' \!\ge 1$ we have 
\begin{align} 
\label{equ:largepmatching}
&\P_p(\cX_{n+1}^y = x, A'_{y,x,n+1}) \le (1-p)n \P_p(\cX_n^y = x),\\
\label{equ:largepnotmatching}
&\P_p(\cX_{n+1}^y = x, (A'_{y,x,n+1})^c) \le c_1^n (1-p)^{\frac n 2 (2m_{y,x}'+1)},
\\
\label{equ:largepprob}
&\P_p(\cX_n^y = x) \ge (1-p)^{c_2} c_3^{-n} (1-p)^{n m_{y,x}'} ,  
\end{align}
where  $A'_{y,x,n+1} := \{\exists 0< i \le \frac {n+1} 2: \cN'_{2i-1} = m'_{y,x} = \cN'_{2i}\}$
\end{Prop}

\Pf This is very similar to the proof of Proposition \ref{prop:smallprep}. Let $b := |E_1|$.
For \eqref{equ:largepmatching} we note that 
by definition of $A'_{y,x,n+1}$ we have 
$$
\P_p(\cX_{n+1}^y = x, A'_{y,x,n+1})
\le \sum_{0 < i \le \frac{n+1}2} \P_p(\cX_{n+1}^y = x, 
\cN'_{2i-1} = m'_{y,x} = \cN'_{2i})
$$
and for every $i$ 
$$
\P_p(\cX_{n+1}^y = x, 
\cN'_{2i-1} = m'_{y,x} = \cN'_{2i})
\le \P_p(\cX_{n}^y = x) (1-p).
$$
To see this, we consider a bond configuration  
$z = (z_e)_{e \in E_{1..n+1}}$ contributing to $\{\cX_{n+1}^{y} = x, \cN'_{2i-1} = m'_{y,x} = \cN'_{2i}\}$. 
We decompose this configuration 
into $(z_e)_{e \in E_{2i-1}}$ and $\bz = (\bz_e)_{e \in E_{1..n}}$, which is obtained from $z$ be deleting layer $E_{2i-1}$. 
By the previous lemma $\bz$ contributes to $\{\cX_n^{y} = x\}$, 
and since $m_{y,x}' \ge 1$ we have 
$$
\P_p(\forall e \in E_{2i-1}: \cZ_e = z_e) 
= (1-p)^{m_{y,x}'} p^{b - m_{y,x}'} \le 1-p.
$$ 
This establishes the above estimate and using 
$\frac{n+1} 2 \le n$ we thus get \eqref{equ:largepmatching}. 
For \eqref{equ:largepnotmatching} we note that 
$$
\P_p(\cX_{n+1}^y = x, (A'_{y,x,n+1})^c)
\le \P_p(\forall 0 < i \le \frac{n+1} 2: \cN'_{2i-1} + \cN'_{2i} > 2 m'_{y,x})  
$$
by definition of $A'_{y,x,n+1}$ and $m'_{y,x}$. 
Furthermore for every $i > 0$ 
\begin{align*}
&\P_p(\cN'_{2i-1} + \cN'_{2i} > 2 m'_{y,x}) 
= \sum_{2 m'_{y,x} + 1 \le k \le 2b} \binom {2b} k (1-p)^kp^{2b-k}\\
&\le \sum_{0 \le k \le 2b} \binom {2b} k (1-p)^{ 2 m'_{y,x} + 1}
= 2^{2b}(1-p)^{ 2 m'_{y,x} + 1}.
\end{align*}
Combining this with the independence of $\cN'_{2i-1} + \cN'_{2i}$ $(i > 0)$, setting $c_1 := 2^{2b}$ and using
$\frac{n+1} 2 \le n, \frac {n+1 } 2 \ge \frac n 2$ 
we obtain \eqref{equ:largepnotmatching}. 
For \eqref{equ:largepprob} we note that by definition of $l'$,  
for some $i \in \{1,...,l'\}$ 
$\{\cX_{l'}^y = x, \cN_i = m'_{y,x}\} \neq \emptyset$. 
Let  $z = (z_e)_{e \in E_{1..l'}}$ be a configuration of bonds 
contributing to this event, and let $\bz = (\bz_e)_{e \in E_{1..n}}$ 
be the bond configuration obtained from $z$ 
by inserting $n-l'$ additional layers above layer $E_i$, 
such that in these additional layers the closed bonds precisely match the positions of the closed bonds in $E_i$. 
By the above lemma $\bz$ contributes to $\{\cX_n^y =x\}$
and thus 
\begin{align*}
&\P_p(\cX_n^y = x) \ge \P_p(\forall e \in E_{1..n}: \cZ_e = \bz_e)\\ 
&= \P_p(\forall e \in E_{1..l'}: \cZ_e = z_e) ((1-p)^{m'_{y,x}}p^{b - m'_{y,x}})^{n-l'}. 
\end{align*}
Noting that $p \ge \frac 1 2$ we estimate 
$\P_p(\forall e \in E_{1..l'}: \cZ_e = z_e) \ge (1-p)^{bl'}$, 
$p^{b-{m'_{y,x}}} \ge (\frac 1 2)^b$
and $n-l' \le n$. Thus we obtain \eqref{equ:largepprob} 
by setting $c_2 := bl'$ and $c_3 := 2^b$.
\qed

If it is possible that all bonds in a layer are open, 
we have to give a separate argument: 
Here with high probability there will be such a layer, 
and we proceed by discarding the layer before that.

\begin{Prop}\label{prop:large2}
There are constants $c_1,c_2,c_3,c_4 \ge 1$ depending only on $G$ such that for all 
$n \ge l'$, $p \in (\frac 1 2,1)$ and 
$y,x \in M^*$ such that $y \to x$ and $m_{y,x}' = 0$ we have 
\begin{align} 
\label{equ:largepfull}
&\P_p(\cX_{n+1}^y = x, B_{n+1}) \le (1-(1-p)^{c_1}) \P_p(\cX_n^y = x),\\
\label{equ:largepnotfull}
&\P_p(\cX_{n+1}^y = x, B_{n+1}^c) \le (c_2(1-p))^n,
\\
\label{equ:largepprobvar}
&\P_p(\cX_n^y = x) \ge c_3^{-n} (1-p)^{c_4} ,  
\end{align}
where $B_{n+1} := \{\exists  2 \le  i \le n+1 \forall e \in E_i: \cZ_e = 1\}$.
\end{Prop}

\Pf Let $b := |E_1|$.
For \eqref{equ:largepfull} let $z = (z_e)_{e \in E_{1..n+1}}$ 
be a bond configuration contributing to $\{\cX_{n+1}^y = x, B_{n+1}\}$. 
Let $i\in \{2,...,n+1\}$ denote the last layer of this configuration 
consisting exclusively of open bonds (which exists by definition of $B_{n+1})$. 
We decompose this configuration into 
$\bz = (\bz_e)_{e \in E_{1..n}}$
which is obtained from $z$ be deleting layer $E_{i-1}$, 
and $z' = (z'_e)_{e \in E_1}$ which equals the bond configuration of the deleted layer (i.e. $z'_e := z_{\tau^{i-2}(e)}$). 
We note that by $z'$ and $\bz$ the original configuration $z$ 
is uniquely determined: If $j \in \{1,...,n\}$ is the last layer of 
$\bz$ consisting of open bonds only, then $z$ can be recovered 
from $\bz$ by inserting $z'$ into $\bz$ between the layers $E_{j-1}$ and $E_j$. 
Furthermore $\bz$ contributes to $\{\cX_n^y = x\}$ (since 
after a layer with all bonds open the pattern is $x_*$ provided the preceding pattern is infected) and $z'$ contributes 
to $\{\exists e \in E_{1}^v: \cZ_e = 1\}$ (since the infection 
has to be transmitted by this layer). 
By independence we thus obtain 
$$
\P_p(\cX_{n+1}^y = x, B_{n+1}) \le \P_p(\cX_n^y = x) \P_p(\exists e \in E_{1}^v: \cZ_e = 1), 
$$
which gives \eqref{equ:largepfull} setting $c_1 := |E_{1}^v| = |V|$. 
For \eqref{equ:largepnotfull} we note that 
$$
\P_p(\cX_{n+1}^y = x, B_{n+1}^c)
\le \P_p(\forall 2 \le i \le n+1: \cN_i \neq b),  
$$
by definition of $B_{n+1}$. 
Furthermore for every $i$ 
$$\P_p(\cN_i \neq  b) = 1 - p^b = (1-p) \sum_{0 \le j < b} p^{j} \le b(1-p). 
$$
Combining this with the independence of $\cN_i$ $(i > 0)$
and setting $c_2 := b$ we obtain \eqref{equ:largepnotfull}. 
For \eqref{equ:largepprobvar} we note that by definition of $l'$,  
for some $i \in \{1,...,l'\}$ 
$\{\cX_{l'}^y = x, \cN_i = b\} \neq \emptyset$. 
Let  $z = (z_e)_{e \in E_{1..l'}}$ be a configuration of bonds 
contributing to this event, and let $(\bz_e)_{e \in E_{1..n}}$ 
be the bond configuration obtained from $z$ 
by inserting $n-l'$ additional layers preceding layer $E_i$, 
such that in these additional layers all bonds are open.
This new configuration contributes to $\{\cX_n^y =x\}$
and thus 
\begin{align*}
&\P_p(\cX_n^y = x) \ge \P_p(\forall e \in E_{1..n}: \cZ_e = \bz_e) 
= \P_p(\forall e \in E_{1..l'}: \cZ_e = z_e) p^{b(n-l')}. 
\end{align*}
Noting that $p \ge \frac 1 2$ we estimate 
$\P_p(\forall e \in E_{1..l'}: \cZ_e = z_e) \ge (1-p)^{bl'}$
and $p^{b(n-l')} \ge (\frac 1 2)^{bn}$, and 
  we obtain \eqref{equ:largepprobvar} by setting $c_3 := 2^b$ and $c_4 = bl'$.
\qed 

\bigskip 

Combining the two preceding propositions gives the desired result in case of large values of $p$: 

\begin{Prop} \label{prop:large}
There is a $p_1 \in (\frac 1 2, 1)$ and an $N \ge 0$ such that 
\begin{equation}
\label{equ:large}
\forall n \ge N \forall p \in [p_1,1) \forall y,x \in M^*: 
\P_p(\cX_n^y =x) \ge \P_p(\cX_{n+1}^y = x). 
\end{equation}
\end{Prop}

\Pf Choosing $c_1,c_2,c_3,c_4 \ge 1$ according to the above proposition and $c_1',c_2',c_3' \ge 1$ according to the previous one, 
we fix $N \ge \max\{l',2(c_4+c_1), 4c_2'\}$ and consider $p > \frac 1 2$ and $y,x \in M^*$. In case of $y \to x$ such that $m_{y,x}' = 0$ 
we get 
\begin{align*}
&\P_p(\cX_{N+1}^y = x) = \P_p(\cX_{N+1}^y = x, B_{N+1}) + \P_p(\cX_{N+1}^y = x, B_{N+1}^c)\\
&\le \P_p(\cX_N^y = x) \big( 1 - (1-p)^{c_1} + (c_2 (1-p))^Nc_3^N(1-p)^{-c_4}\big) \\
&\le \P_p(\cX_N^y = x)
\big(1 - (1-p)^{c_1}(1- (c_2c_3)^N (1-p)^{\frac N 2})\big). 
\end{align*}
If $1-p$ is sufficiently small, we thus have 
$\P_p(\cX_{N+1}^y = x) \le \P_p(\cX_N^y = x)$. 
In case of $y \to x$ such that $m_{y,x}' \ge 1$ we get 
\begin{align*}
&\P_p(\cX_{N+1}^y = x) = \P_p(\cX_{N+1}^y = x, A'_{y,x,N+1}) + \P_p(\cX_{N+1}^y = x, (A'_{y,x,N+1})^c)\\
&\le \P_p(\cX_N^y = x) \big( (1-p) N + (c_1')^N (1-p)^{\frac N 2} (c_3')^{N} (1-p)^{-c_2'}\big)\\
&\le \P_p(\cX_N^y = x)
\big( (1-p) N + (c_1'c_3')^N (1-p)^{\frac N 4}\big). 
\end{align*}
If $1-p$ is sufficiently small, we thus have 
$\P_p(\cX_{N+1}^y = x) \le \P_p(\cX_N^y = x)$. 
In case of $y \not \to x$ we have 
$\P_p(\cX_{N+1}^y = x) = 0 =  \P_p(\cX_N^y = x)$. 
In conclusion there is a value $p_1 \in (\frac 1 2,1)$ 
(depending on $G$ only) such that for all $p \in [p_1,1)$ and 
$y,x \in M^*$ we have 
$\P_p(\cX_N^y = x) \ge \P_p(\cX_{N+1}^y = x)$, and 
\eqref{equ:large} follows using Lemma \ref{lem:monoN}. 
\qed

\subsection{Proof of Theorem \ref{Thm:pattern}}

Theorem \ref{Thm:pattern} can now easily be proved combining 
the results from the last three subsections: 
We first choose $p_0 \in (0,\frac 1 2)$ and $N_1 \ge 0$ 
according to Proposition \ref{prop:small}, 
and $p_1 \in (\frac 1 2, 1)$ and $N_2 \ge 0$ 
according to Proposition \ref{prop:large}. 
For these values of $p_0,p_1$ we choose $N_3 \ge 0$ 
according to Proposition \ref{prop:intermediate}. 
Setting $N := \max \{N_1,N_2,N_3\}$, 
for all $n \ge N, p \in (0,1)$, $y,x \in M^*$ we thus get 
$\P_p(\cX_n^y = x)$ $\ge \P_p(\cX_{n+1}^y = x)$ as desired. 
The two additional claims easily follow from the first one 
using Proposition \ref{Prop:monodiff}.
\qed

\section{Number of infected points per layer} \label{Sec:Number}

Here we prove Theorem \ref{Thm:number}. 
Let $G$ be a connected graph of bounded degree such that $\deg(x) \le \De$ for all $x \in  V$. 
We first show that for $p \in [0,1]$ and $n \ge 0$ we have 
\begin{align} 
\label{equ:numbermono1}
\E_p(\cW_n) \cdot \max_{w \in G} \E_p(\tcW^w_1) \ge \E_p(\cW_{n+1}),   
\end{align}
where $\tcW^w_1 = |\{v \in V\!: (v,1) \lra_{\ge 1} (w,0)\}|$
and $\lra_{\ge n}$ is the connectivity relation  induced 
by percolation on $E_{n..}$.  
For this we note that 
$$
\E_p(\cW_{n+1}) = \sum_{v \in V} \P_p(o' \lra (v,n+1)),
$$
and considering a self-avoiding path from $o'$ to $(v,n+1)$ and its last vertex in the layer $V_n$, 
we see that
$$
\{o' \lra (v,n+1)\} =  \bigcup_{w \in V} \{o' \lra (w,n)\}  \circ 
\{(w,n) \lra_{\ge n+1} (v,n+1)\}.
$$
Thus using a union bound and the BK-inequality (see \cite{G}) 
we obtain 
\begin{align*}
\E_p(\cW_{n+1}) &\le \sum_{w,v \in V} \P_p(o' \lra (w,n)) 
\P_p((w,n) \lra_{\ge n+1} (v,n+1)).
\end{align*}
While the BK-inequality is usually formulated 
in a setting with finitely many bonds, 
this is not a problem since the above estimate thus
holds for percolation on $E_{-n'..n+1+n'}$ 
for any fixed $n' \ge 0$, and letting $n' \to \infty$ all probabilities and expectations converge to give the above estimate as stated. 
By translation invariance 
\begin{align*}
&\sum_{v \in V} \P_p((w,n) \lra_{\ge n+1} (v,n+1))
= \sum_{v \in V} \P_p((w,0) \lra_{\ge 1} (v,1)) = \E_p(\tcW_1^w)
\end{align*}
and we thus obtain \eqref{equ:numbermono1}. 
\eqref{equ:numbermono} follows from \eqref{equ:numbermono1} by estimating $\E_p(\tcW_1^w)$. Let $w \in V$ be fixed for the remainder of the proof. 
Let $W_{l}$ denote the set of self-avoiding paths in $E_{1..}$ of length $l$ 
from $(w,0)$ to some vertex in $V_1$, 
and for any $P \in W_{l}$ let $A_P$ denote the event, that all edges of $P$ are open. 
We note that 
$$
\E_p(\tcW_1^w) \le \E_p(\sum_{l \ge 1} \sum_{P \in W_{l}} 1_{A_P}) = \sum_{l \ge 1} \sum_{P \in W_{l}}
\P_p(A_P) = \sum_{l \ge 1} |W_{l}| p^l.
$$
Any $P \in W_{l}$ starts with a step going up (to $V_1$) and then never returns to $V_0$ (because it is self-avoiding). 
It is possible, that it stays within $V_1$, in which case 
we have  $\De$ possibilities for the first step and $\De -1$ 
possibilities for all other steps. 
On the other hand, if the path also visits $V_2$, 
then (after the first step up) it remains in $V_1$ 
for $k_1 \ge 0$ steps, then goes up to $V_2$, then performs 
another $k \ge 1$ steps to end in $V_2$, then goes down to $V_1$, 
then stays in $V_1$ for an additional $k_2 \ge 0$ steps. 
We note that for the last step of the $k$ steps from $V_2$ to $V_2$ 
there are at most $\De$ choices.  (After $k-1$ steps we are either in $V_1,V_2$ or $V_3$. If we are in $V_1$ or $V_3$ we have to 
go up or down respectively. If we are in $V_2$ we stay in $V_2$
by $\De$ choices.) This gives 
$$
|W_{l}| \le \De(\De-1)^{l-2}  + \sum_{k_1,k,k_2}
 \De(\De-1)^{k_1-1} \cdot (\De+1)^{k-1}\De \cdot \De(\De-1)^{k_2-1}, 
$$
where the sum is over all $k_1 \ge 0, k \ge 1, k_2 \ge 0$ 
such that $k_1+k+k_2+3 = l$, and we interpret 
$\De (\De-1)^{-1}$ as $1$. We thus obtain 
\begin{align*}
\E_p(\tcW_1^w) &\le \sum_{l \ge 1} \De(\De-1)^{l-2} p^{l}  \\&
+ \sum_{k_1 \ge 0} \De(\De-1)^{k_1-1} p^{k_1+ 1}
 \sum_{k \ge 1} (\De+1)^{k-1}\De p^{k+1} 
  \sum_{k_2 \ge 0} \De(\De-1)^{k_2-1} p^{k_2+1} \\
&= \frac{p(1+p)}{1 - (\De - 1)p} + (\frac{p(1+p)}{1 - (\De - 1)p})^2 
\frac{\De  p^2}{1 - (\De +1)p} \quad  =: \quad f_\De(p) 
\end{align*}
using $p \le  p_\De := \frac 1 {\De + 1.4} < \frac 1 {\De + 1}$. 
Since $f_\De(p)$ is increasing in $p$, 
it suffices to check that $g(\De) := f_\De(p_\De) \le 1$. 
Using $\frac 1 {p_\De} - (\De - 1) = 2.4$ and $\frac 1 {p_\De} - (\De + 1) = 0.4$ we have 
\begin{align*}
&g(\De) = \frac{1+\frac 1 {\De + 1.4}}{2.4} + (\frac{1+\frac 1 {\De + 1.4}}{2.4})^2 \frac{\frac \De {\De + 1.4}}{0.4}
\le \frac{1+\frac 1 {\De + 1.4}}{2.4} (1 + \frac 1 {2.4\cdot 0.4}) =: h(\De) 
\end{align*}
using $(1 + \frac 1 {\De + 1.4})\frac{\De}{\De + 1.4} \le 1$. 
$h$ is decreasing in $\De$, and it is easy to check 
that $g(\De) \le 1$ for $\De \in \{0,1,2,3,4\}$ and $h(5) \le 1$ , thus $g(\De) \le 1$ for all $\De$ and we are done. 
\qed

\end{document}